\documentclass[12pt]{amsart}
\usepackage{graphicx}
\usepackage{amsfonts}
\usepackage{amssymb}
\usepackage{latexsym}
\usepackage{amscd}

\renewcommand{\marginpar}[1]{}
\catcode`\@=12

\def\Empty{}
\newcommand\oplabel[1]{
  \def\OpArg{#1} \ifx \OpArg\Empty {} \else
  	\label{#1}
  \fi}
		
\newcommand{\comm}[1]{}

\input{psfig}
\newtheorem{thm}{Theorem}[section]
\newtheorem{cor}[thm]{Corollary}
\newtheorem{lem}[thm]{Lemma}
\newtheorem{prop}[thm]{Proposition}

\newenvironment{pf*}[1]{\proof[#1]}{\endproof}
\usepackage{euscript}

\usepackage[OT2,OT1]{fontenc}

\newcommand{\cal}[1]{{\mathcal #1}}

\newcommand{\beq}{\begin{equation}}
\newcommand{\eeq}{\end{equation}}
\newcommand{\eref}[1]{(\ref{#1})}

\newcommand{\ve}{\varepsilon}
\newcommand{\de}{\delta}
\newcommand{\al}{\alpha}
\newcommand{\be}{\beta}
\newcommand{\ga}{\gamma}

\newcommand{\om}{\omega}

\theoremstyle{definition}
\newtheorem{defn}{Definition}[section]

\theoremstyle{remark}

\newcommand{\dist}{\operatorname{dist}}

\newcommand{\eps}{\epsilon}

\numberwithin{equation}{section}
\newcommand{\thmref}[1]{Theorem~\ref{#1}}

\newcommand{\secref}[1]{\S\ref{#1}}
\newcommand{\lemref}[1]{Lemma~\ref{#1}}

\newcommand{\cF}{{\cal F}}

\newcommand{\cB}{{\cal B}}

\newcommand{\cC}{{\cal C}}

\newcommand{\CC}{{\mathbb C}}
\newcommand{\RR}{{\mathbb R}}
\newcommand{\TT}{{\mathbb T}}

\newcommand{\NN}{{\mathbb N}}
\newcommand{\DD}{{\mathbb D}}

\newcommand{\QQ}{{\mathbb Q}}

\begin{document}
\addtolength{\evensidemargin}{-0.7in}
\addtolength{\oddsidemargin}{-0.7in}

\title[Complexity of Julia sets]{On computational complexity of Siegel Julia sets}
\author{I. Binder, M. Braverman, M. Yampolsky}
\thanks{The first and third authors are partially supported by NSERC Discovery grants.
The second author is partially supported by NSERC Postgraduate Scholarship}

\date{\today}
\begin{abstract}
It has been previously shown by two of the authors
 that some polynomial Julia sets are algorithmically impossible to draw with 
arbitrary magnification. On the other hand, for a large class of examples the problem
of drawing a picture has polynomial complexity. In this paper we demonstrate the 
existence of computable quadratic Julia sets whose computational complexity is arbitrarily high.
\end{abstract}
\maketitle

%\begin{abstract}
%\end{abstract}
\section{Foreword}
Let us informally say that a compact set in the plane is computable if one
can program a computer to draw a picture of this set on the screen,
with an arbitrary desired magnification. It was recently shown by the 
second and third authors, that some Julia sets are not computable \cite{BY}.
This in itself is quite surprising to dynamicists -- Julia sets are among the
``most drawn'' objects in contemporary mathematics, and numerous algorithms
exist to produce their pictures. 
In the cases when one has not been able to produce informative pictures
(the dynamically pathological cases, like maps with a Cremer or a highly
Liouville Siegel point) the feeling had been that this was due to
the immense computational resources required by the known algorithms.

The next surprise came with the discovery by the authors of this paper in
\cite{BBY} that all Cremer quadratics (or more generally, rational maps
without rotation domains) have computable Julia sets.
The non-computable examples
constructed in \cite{BY} were Siegel quadratic polynomials, and one 
would expect the Cremer case to be at least as bad if not worse 
computationally.

The natural question to ask is then whether in those  cases in which we know
the Julia set is computable, but no good pictures exist, the computational
complexity of such a set is indeed high. Here at least, our original intuition
seems to be correct: it is shown in the present paper that there exist
computable Siegel quadratic Julia sets with arbitrarily high computational 
complexity. An irritating possibility still remains that some Cremer Julia
sets are computationally easy (and we just do not go about trying to draw
them in the right way). This, however, seems unlikely.
We note that the  examples constructed in this paper are the first known cases
of Julia sets which are not poly-time computable. 
The second author \cite{thesis} and independently
Rettinger \cite{Ret} have previously shown that hyperbolic Julia sets are
poly-time computable. More recently the second author has shown \cite{Brv2} that some
 Julia sets with parabolics are poly-time computable as well. The last result
was yet another surprise, as the time complexity of all previously known 
algorithms for these Julia sets was exponential.

The structure of the paper is as follows. In \secref{sec:main thm} of the Introduction,
having stated the principal definitions, we formulate the main result of the paper.
In \secref{sec:outline} we give a sketch of the argument. In \secref{sec changes} we
prove several technical lemmas. The final \secref{sec:proof}
contains the proof of the Main Theorem.

\subsection*{Acknowledgement.} We would like to thank Giovanni Gallavotti for 
very helpful suggestions on the exposition.

\section{Introduction}
\label{section-intro}

\subsection{Computability of real sets}
The reader is directed to \cite{BY} for a more detailed
 discussion of the notion of computability of
subsets of $\RR^n$ as applied, in particular, to Julia sets. We recall the 
principal definitions here. The exposition below uses the concept of a {\it Turing
Machine}. This is a standard model for a computer program employed by computer
scientists. Readers unfamiliar with this concept should think instead of an
algorithm written in their favorite programming language. These concepts are
known to be equivalent.

Denote by $\DD$ the set of the {\it dyadic rationals}, that is, rationals of 
the form $\frac{p}{2^m}$. 
%Rationals in $\DD$ can be easily represented as binary numbers. 
We say that $\phi: \NN \rightarrow \DD$ is an {\it oracle}
for a real number $x$, if $| x - \phi(n)|<2^{-n}$ for all $n \in \NN$. 
In other words, $\phi$ provides a good dyadic approximation for $x$. 
We say that a Turing Machine (further abbreviated as TM) $M^{\phi}$ is 
an {\it oracle machine}, if at every step of 
the computation M is allowed to query the value $\phi(n)$ for any $n$.
This definition allows us to define the computability of real functions 
on compact sets. 

\begin{defn}
\label{funcomp}
We say that a function $f:[a,b] \rightarrow [c,d]$ is computable, if 
there exists an oracle TM $M^{\phi} (m)$ such that if $\phi$ is an oracle 
for $x \in [a,b]$, then on input $m$, $M^{\phi}$ outputs a $y \in \DD$ 
such that $| y - f(x)|<2^{-m}$.
\end{defn}

\noindent
To understand this definition better, the reader without a Computer Science
background should think of a computer program with an instruction
$$\text{ READ real number }x\text{ WITH PRECISION }n(m).$$
On the execution of this command, a dyadic rational $d$ is input from the 
keyboard. This number must not differ from $x$ by more than  $2^{-n(m)}$ 
(but otherwise can be arbitrary). The algorithm then outputs $f(x)$ to precision
$2^{-n}$.

\noindent
It is worthwhile to note why the oracle mechanism is introduced. 
There are only countably many possible algorithms, and consequently only
countably many {\it computable} real numbers which such algorithms can
encode. Therefore, one wants to separate the hardness of encoding the
real number $x$ from the hardness of computing the value of the function
$f(x)$, having the access to the value of $x$.

\noindent
Let $K \subset \RR^k$ be a compact set. 
We say that a TM M computes the set $K$ if it approximates $K$ in the {\it 
Hausdorff metric}. Recall that the Hausdorff metric is a metric on 
compact subsets of $\RR^n$ defined by 
\begin{equation}
\label{hausdorff metric}
d_H ( X, Y) =  \inf \{\epsilon > 0 \;|\; X \subset U_{\epsilon} 
(Y)~~\mbox{and}~~  Y \subset U_{\epsilon}(X)\},
\end{equation}

\noindent
where $U_\eps(S)$ is defined as the union of the set of $\eps$-balls with centers in $S$.

\noindent
 We introduce
a class $\cC$ of sets which is dense in metric $d_H$ among the compact 
sets and which has a natural correspondence to binary strings. 
Namely $\cC$ is the set of finite unions of dyadic balls:
$$
\cC= \left\{ \bigcup_{i=1}^n \overline{B(d_i, r_i)}~|~~\mbox{where}~~d_i, 
r_i 
\in \DD \right\}.
$$
Members of $\cC$ can be encoded as binary strings in a natural way. 

\noindent
We now define the notion of computability of subsets of $\RR^n$ (see \cite{Wei}, and
 also \cite{WeiPaper}).

\begin{defn}
\label{setcomp}
We say that a compact set $K\subset \RR^k$ is computable, if there exists a
TM $M(d,n)$, where $d\in\DD$, $n\in \NN$ which outputs a value $1$ if $\dist(d,K)<2^{-n}$,
the value $0$ if $\dist(d,K)>2\cdot 2^{-n}$, and in the ``in-between'' case it halts and
outputs either $0$ or $1$. 

In other words, it
computes, in the classical sense,  a function from the family $\cF_K$
of functions of the form
\begin{equation}
\label{defcomp}
f(d,n)=\left\{
\begin{array}{ll}
0,& \text{if }\dist(d,K)>2\cdot 2^{-n}\\
1,& \text{if }\dist(d,K)<2^{-n}\\
0\text{ or }1,&\text{otherwise}
\end{array}
\right.
\end{equation} 
\end{defn}

\begin{thm}
\label{equiv-defn}
For a compact $K \subset \RR^k$ the following are equivalent:

(1) $K$ is computable as per definition \ref{setcomp},

(2) there exists a TM 
$M(m)$, such that on input $m$, $M(m)$ outputs an encoding of $C_m \in 
\cC$ such that $d_H (K, C_m) < 2^{-m}$ (global computability),

(3) the {\em distance function} $d_K (x) = \inf \{ |x-y|~~|~~y\in K \}$ 
is 
computable as per definition \ref{funcomp}.
\end{thm}

Note that in the case $k=2$  computability means that $K$ can be drawn on a computer screen 
with arbitrarily good precision (if we imagine the screen as a lattice of pixels).

In the present paper we are interested in   questions concerning 
the computability of the Julia set $J_c = J(f_c) = J(z^2+c)$.
Since there 
are uncountably many possible parameter values for $c$, we cannot 
expect for each $c$ to have a machine $M$ such that $M$ computes 
$J_c$ (recall that there are countably many TMs). On the other hand, 
it is reasonable to want $M$ to compute $J_c$ with an oracle 
access to $c$. Define the function $J: \CC \rightarrow K^{*}$ ($K^{*}$ is 
the set of all compact subsets of $\CC$) by $J(c)=J(f_c)$. 
In a complete analogy to Definition \ref{funcomp} we can define

\begin{defn}
\label{funcomp2}
We say that a function $\kappa:S \rightarrow K^{*}$ for some bounded set $S$ is 
computable, if there exits an oracle TM $M^{\phi}(d,n)$, 
where $\phi$ is an oracle for $x\in S$, which computes a function (\ref{defcomp}) of the family
$\cF_{\kappa(x)}$.

Equivalently, there exists an oracle TM $M^\phi(m)$ with $\phi$ again representing
$x\in S$ such that 
 on input $m$, $M^{\phi}$ outputs a $C \in \cC$ 
such that $d_H (C, \kappa(x))<2^{-m}$.
\end{defn}

\noindent
In the case of Julia sets:

\begin{defn}
\label{Jcomp}
We say that $J_c$ is computable if the function $J: d \mapsto J_d$ is 
computable on the set $\{ c \}$.
\end{defn}

\noindent
We have the following (see \cite{thesis}):
\begin{thm}
Suppose that a TM $M^\phi$ computes the function $J$ on a set $S\subset\CC$. Then $J$ is continuous
on $S$ in Hausdorff sense.
\end{thm}

\begin{proof}
Let $c$ be any point in $S$, and let $\ve = 2^{-k}$ be given. Let $\phi$ be 
an oracle for $c$ such that $|\phi(n)-c|<2^{-(n+1)}$ for all $k$. We run 
$M^\phi (k+1)$ with this oracle $\phi$. By the definition of $J$, it outputs a set $L$ which
is a $2^{-(k+1)}$ 
approximation of $J_c$ in the Hausdorff metric. 

The computation is performed in a finite amount of time. Hence 
there is an $m$ such that $\phi$ is only queried with parameters not exceeding $m$. 
Then for any $x$ such that $|x-c|<2^{-(m+1)}$, $\phi$ is a valid oracle for $x$ 
up to parameter value of $m$. In particular, we can create an oracle $\psi$ for $x$ 
that agrees with $\phi$ on $1, 2, \ldots, m$. If $x \in S$, then the execution 
of $M^{\psi}(k+1)$ will be identical to the execution of $M^{\phi}(k+1)$, and it 
will output $L$ which has to be an approximation of $J_x$. Thus we have
$$
d_H(J_c, J_x) \le d_H (J_c, L)+d_H (J_x, L) < 2^{-(k+1)}+2^{-(k+1)} = 2^{-k}.
$$
This is true for any $x \in B(c,2^{-(m+1)})\cap S$. Hence $J$ is continuous
on $S$. 
\end{proof}

\noindent
The second and third authors have demonstrated in \cite{BY}:

\begin{thm}
\label{non-comput}
There exists a parameter value $c\in \CC$ such that the Julia set of the
quadratic polynomial $f_c(z)=z^2+c$ is not computable.
\end{thm}

\noindent
The quadratic polynomials in \thmref{non-comput} possess Siegel disks (see \secref{sieg-quad}
below for the definitions of Siegel and Cremer points).
It was further shown by the authors of the present paper in \cite{BBY}
that the absence of rotation domains, that is either Siegel disks or
Herman rings, guarantees computability of the rational Julia set.
This implies, in particular, that all Cremer quadratic Julia sets are computable --
this despite the fact that no informative high resolution images of such
sets have ever been produced. One expects, however, that such ``bad'' but
still computable examples have high algorithmic complexity, which makes
the computational cost of producing such a picture prohibitively high. We note,
that the second author \cite{thesis} and independently Rettinger \cite{Ret}
have shown:

\begin{thm}
Hyperbolic Julia sets are computable in polynomial time. That is, if $J$ is
the Julia set of a hyperbolic rational mapping $R$, then there exists a TM
$M(d,n)$  which computes a function
of the family  (\ref{defcomp}) in time polynomial in the bit size of the input $(d,n)$.
It is worth noting that the same oracle TM $M^\phi(d,n)$
with the oracle representing the parameters of the rational mapping $R$,
can be selected for all hyperbolic Julia sets of the same degree. 
Moreover, the asymptotics of the
polynomial time bound  depends only on $R$ but not on the input $(d,n)$.
\end{thm}

\subsection{Statement of the Main Theorem}
\label{sec:main thm}

\noindent
On the other end of the complexity spectrum we expect to find ``bad'' but computable Siegel Julia
sets and Cremer Julia sets. Indeed, it the present paper we show:

\begin{thm}
There exist quadratic Siegel Julia sets of arbitrarily high computational complexity.
More precisely, for any computable increasing function $h:\NN\to\NN$ there exists a computable
Siegel parameter value $c\in\CC$ such that:
\begin{itemize}
\item the Julia set $J_c$ is computable by an oracle TM;
\item for any oracle TM $M^\phi(m)$ which computes the $2^{-m}$-approximations
to  $J_c$, 
there exists a sequence $\{m_i\}_{i=1}^\infty$ such that
$M^\phi$  requires the time of at least $h(m_i)$ to compute the 
approximation $C_{m_i}\in\cC$.
\end{itemize} 
\end{thm}

\noindent
From this statement for global computational complexity immediately follows the corresponding local
statement:

\begin{cor}
There exist computable parameter values $c$ for which the Julia set $J_c$ is computable,
and the complexity of the problem of computing a function (\ref{defcomp}) in the family $\cF_{J_c}$
is arbitrarily high.
\end{cor}

\subsection{Siegel disks of quadratic maps}
\label{sieg-quad}
Let $R:\hat\CC\to\hat\CC$ be a rational map of the Riemann sphere.
For a periodic point $z_0=R^p(z_0)$
of period $p$ its {\it multiplier} is the quantity $\lambda=\lambda(z_0)=DR^p(z_0)$.
We may speak of the multiplier of a periodic cycle, as it is the same for all points
in the cycle by the Chain Rule. In the case when $|\lambda|\neq 1$, the dynamics
in a sufficiently small neighborhood of the cycle is governed by the Intermediate 
Value Theorem: when $0<|\lambda|<1$, the cycle is {\it attracting} ({\it super-attracting}
if $\lambda=0$), if $|\lambda|>1$ it is {\it repelling}.
Both in the attracting and repelling cases, the dynamics can be locally linearized:
\begin{equation}
\label{linearization-equation}
\psi(R^p(z))=\lambda\cdot\psi(z)
\end{equation}
where $\psi$ is a conformal mapping of a small neighborhood of $z_0$ to a disk around $0$.
By a classical result of Fatou, a rational mapping has at most finitely many non-repelling
periodic orbits.

In the case when $\lambda=e^{2\pi i\theta}$, $\theta\in\RR$, 
 the simplest to study is the {\it parabolic case} when $\theta=n/m\in\QQ$, so $\lambda$ 
is a root of unity. In this case $R^p$ is not locally linearizable; it is not hard to see that $z_0\in J(R)$.
 In the complementary situation, two non-vacuous possibilities  are considered:
{\it Cremer case}, when $R^p$ is not linearizable, and {\it Siegel case}, when it is.
In the latter case, the linearizing map $\psi$ from (\ref{linearization-equation}) conjugates
the dynamics of $R^p$ on a neighborhood $U(z_0)$ to the irrational rotation by angle $\theta$
(the {\it rotation angle})
on a disk around the origin. The maximal such neighborhood of $z_0$ is called a {\it Siegel disk}.

Let us discuss in more detail the occurrence of Siegel disks in the quadratic family.
For a number $\theta\in [0,1)$ denote $[r_0,r_1,\ldots,r_n,\ldots]$, $r_i\in\NN\cup\{\infty\}$ its possibly finite 
continued fraction expansion:
\begin{equation}
\label{cfrac}
[r_0,r_1,\ldots,r_n,\ldots]\equiv\cfrac{1}{r_0+\cfrac{1}{r_1+\cfrac{1}{\cdots+\cfrac{1}{r_n+\cdots}}}}
\end{equation}
Such an expansion is defined uniquely if and only if $\theta\notin\QQ$. In this case, the {\it rational 
convergents } $p_n/q_n=[r_0,\ldots,r_{n-1}]$ are the closest rational approximants of $\theta$ among the
numbers with denominators not exceeding $q_n$. In fact, setting $\lambda=e^{2\pi i\theta}$, we have
$$|\lambda^h-1|>|\lambda^{q_n}-1|\text{ for all }0<h<q_{n+1},\; h\neq q_n.$$
The difference $|\lambda^{q_n}-1|$ lies between $2/q_{n+1}$ and $2\pi/q_{n+1}$,
therefore the rate of growth of the denominators $q_n$ describes how well 
$\theta$ may be approximated with rationals.

We recall a theorem due to Brjuno (1972):
\begin{thm}[\cite{Bru}]
Let $R$ be an analytic map with a periodic point $z_0\in\hat\CC$. Suppose that the 
multiplier of $z_0$ is $\lambda=e^{2\pi i\theta}$, and
\begin{equation}
\label{brjuno}
B(\theta)=\displaystyle\sum_n\frac{\log(q_{n+1})}{q_n}<\infty.
\end{equation}
Then $z_0$ is a Siegel point.
\end{thm}

\noindent
Note that a quadratic polynomial with a fixed Siegel disk with rotation angle $\theta$ after an affine
change of coordinates can be written as 
\begin{equation}
\label{form-1}
P_\theta(z)=z^2+e^{2\pi i \theta}z.
\end{equation}
\noindent
In 1987 Yoccoz \cite{Yoc} proved the following converse to Brjuno's Theorem:

\begin{thm}[\cite{Yoc}]
Suppose that for $\theta\in[0,1)$ the polynomial $P_\theta$ has a Siegel point at the origin.
Then $B(\theta)<\infty$.
\end{thm}

\noindent
The numbers satisfying (\ref{brjuno}) are called Brjuno numbers; the set of all Brjuno numbers will be denoted $\cB$. It is a full measure set which contains all Diophantine rotation numbers.
In particular, the rotation numbers $[r_0,r_1,\ldots]$ of {\it bounded type},
that is with $\sup r_i<\infty$ are in $\cB$.
The sum of the series (\ref{brjuno}) is called the Brjuno function. 
For us a different characterization of $\cB$ will be more useful. Inductively define $\theta_1=\theta$
and $\theta_{n+1}=\{1/\theta_n\}$. In this way, 
$$\theta_n=[r_{n-1},r_n,r_{n+1},\ldots].$$
We define the {\it Yoccoz's Brjuno function} as
$$\Phi(\theta)=\displaystyle\sum_{n=1}^{\infty}\theta_1\theta_2\cdots\theta_{n-1}\log\frac{1}{\theta_n}.$$
One can verify that $$B(\theta)<\infty\Leftrightarrow \Phi(\theta)<\infty.$$
The value of the function $\Phi$ is related to the size of the Siegel disk in the following way.

\begin{defn}
Let $(U,u)$ be a simply-connected subdomain of $\CC$ with a marked interior point.
Consider the unique conformal isomorphism $\phi:\DD\mapsto U$ with
$\phi(0)=u$, and $\phi'(0)>0$. The {\it conformal radius of $(U,u)$} is the value of the
derivative $r(U,u)=\phi'(0)$.

Let $P(\theta)$ be a quadratic polynomial with a Siegel disk $\Delta_\theta\ni 0$. 
The {\it conformal radius of the Siegel disk $\Delta_\theta$} is
$r(\theta)=r(\Delta_\theta,0)$.
For all other $\theta\in[0,\infty)$ we set $r(\theta)=0$, and $\Delta_\theta=\{0\}$. 
\end{defn} 

\noindent
By the Koebe 1/4 Theorem of classical complex analysis (see e.g. \cite{Ahlfors}), the  radius of 
the largest Euclidean disk around $u$ which can be inscribed in $U$ is at least $r(U,u)/4$.

\noindent
We note that one has the following direct consequence of the Carath{\'e}odory Kernel Theorem
(see e.g. \cite{Pom}):

\begin{prop}
\label{radius-continuous}
The conformal radius of a quadratic Siegel disk varies continuously with respect to the Hausdorff 
distance on Julia sets.
\end{prop}

\comm{ \noindent
Moreover, from ??? it is not hard to show that the conformal radius of the Siegel 
disk can be computed from the Julia set in bounded time.

\begin{prop}
\label{radius-computable}
Suppose that a Siegel Julia set $J$ is computable in time $T(n)$. Then
the conformal radius of the Siegel disk is (uniformly) computable in
time $T(2^{2^n})$???
\end{prop}}

\noindent
Yoccoz \cite{Yoc} has shown that the sum 
$$\Phi(\theta)+\log r(\theta)$$
is bounded below independently of $\theta\in\cB$. Recently, Buff and Ch{\'e}ritat have greatly improved this result
by showing that:

\begin{thm}[\cite{BC}]
\label{phi-cont}
The function $\theta\mapsto \Phi(\theta)+\log r(\theta)$ extends to $\RR$ as a 1-periodic continuous
function.
\end{thm}

\medskip
\noindent
In \cite{BBY} we obtain the following result on 
computability of quadratic Siegel disks:

\begin{thm}
\label{thm BBY}
The following statements are equivalent:
\begin{itemize}
\item[(I)] the Julia set $J(P_\theta)$ is computable;
\item[(II)] the conformal radius $r(\theta)$ is computable;
\item [(III)] the inner radius $\inf_{z\in\partial \Delta_\theta}|z|$ is computable.
\end{itemize}
\end{thm}

\noindent
We note that when $\theta$ is not a Brjuno number, the quantities in (II) and (III) are each
equal to zero, and the claim is simply that $J(P_\theta)$ is computable in this case.

\medskip

\noindent
We will make use of the following Lemma which bounds
the variation of the conformal radius under a perturbation of the domain.
It  is a direct consequence of the Koebe Theorem (see e.g. \cite{RZ} for a
proof).

\begin{lem}
\label{rad-modulus1}
Let $U$ be a simply-connected subdomain of $\CC$ containing the point $0$ in the
interior. Let $V \subset U$ be a subdomain of $U$. 
Assume that $\partial V\subset B_\eps(\partial U)$.
Then 
$$0<r(U,0)-r(V,0)\leq 4 \sqrt{r(U, 0)}\sqrt{\eps}.$$
\end{lem}

\subsection{Outline of the construction.}
\label{sec:outline}
We can now describe the idea of our construction.
This outline is rather sketchy and suffers from obvious logical deficiencies,
however, it presents the construction in a simple to understand form.
Consider the oracle Turing machines $M^\phi$ with $\phi$ representing the 
parameter $\theta$ in $P_\theta$.
Since there are only countably many Turing machines, we may order these machines in a 
sequence $M^{\phi}_1,$ $M^\phi_2,\ldots$ We denote $S_i$ the domain on which 
 $M^{\phi}_i$ computes $J({P_\theta})$ 
properly. 
We thus have that 
for each $i$, the function $J: \theta \mapsto J({P_\theta})$ is continuous on 
$S_i$.

\noindent
Let us start with a machine $M^\phi_{n_1}$ which computes $J(P_{\theta_*})$ for $\theta_*=[1,1,1,\ldots]$.
If any of the digits $r_i$ in this infinite continued fraction is changed to a sufficiently large $N\in\NN$,
the conformal radius of the Siegel disk will become small. For $N\to\infty$ the Siegel disk will implode
and its center will become a parabolic fixed  point in the Julia set (see \cite{Do}).

If we are careful, we may select $i_1>1$ and $N_1>>1$ in such a way, that for $\theta_1$
given by the continued fraction
where all digits are ones except $r_{i_1}=N_1$ we have
\beq \label{exam1} r(\theta_*)(1-1/4)<r(\theta_1)<r(\theta_*)(1-1/8).\eeq
By the Koebe 1/4-Theorem, there exists $\ell_1>0$ such that the distance between the two Julia sets 
$$d_H(J(P_{\theta_*}),J(P_{\theta_1}))>2^{-\ell_1}.$$
To ensure that the machine $M_{n_1}^\phi$ will not be able to produce an accurate $2^{-\ell_1}$-approximation
of $J(P_{\theta_1})$ faster than in the time $h(\ell_1)$ we simply select $i_1>h(\ell_1)$. This guarantees
that the TM will have to read at least $h(\ell_1)$ digits of the oracle $\phi$ to distinguish the two Julia
sets, which takes the time $h(\ell_1)$.

To ``fool'' the machine $M_{n_2}^\phi$ we then change a digit $r_{i_2}$ for $i_2>i_1$ sufficiently far in the 
continued fraction of $\theta_1$ to a large $N_2$. In this way, we will obtain a Brjuno number $\theta_2$
for which
\beq \label{exam2} r(\theta_*)(1-1/4-1/8)<r(\theta_2)<r(\theta_*)(1-1/4).\eeq
Again, there exists $\ell_2$ such that for any such Brjuno number, we have 
$$d_H(J(P_{\theta_1}),J(P_{\theta_2}))>2^{-\ell_2},$$
and we choose $i_2>h(\ell_2)$. Continuing inductively, we arrive at the desired limiting Brjuno number $\theta_\infty$.
 
To convince the reader that this construction is not artificial, and not due to the peculiarities of the
selected computation model let us recast it somewhat informally as follows. It is possible by an arbitrarily
small perturbation of the parameter $\theta$ to cause a detectable disturbance in the picture of $J(P_\theta)$.
To distinguish the picture of the new Julia set from the old one, in practice one needs to draw it with {\it arbitrary precision
arithmetic}. That is, not only the input of the parameter (reading the oracle) will take a long time due to the number of significant digits, but also all the arithmetic manipulations with this parameter. Of course, 
the former consideration is already sufficient
to prove the theorem.

\section{Computing Noble Siegel Disks}

The primary goal of the present paper is to show that there 
are computationally hard yet {\em computable} Julia sets 
with Siegel disks. To establish this computability 
we need a computability result for {\it noble} Siegel disks.
The term ``noble'' is applied in the literature to rotation numbers
of the form $[a_0,a_1,\ldots,a_k,1,1,1,\ldots]$. The noblest of 
all is the golden mean $\gamma_*=[1,1,1,\ldots]$.

\begin{lem}
\label{bounded-type}
There is a Turing Machine $M$, which given a finite sequence
of numbers $[a_0, a_1, \ldots, a_k]$ computes the conformal radius
$r_\gamma$ for 
the noble number $\gamma=[a_0,  \ldots, a_k, 1,  \ldots]$. 
\end{lem}

\medskip
\noindent
%We recall first how in \cite{BBY} we have established computability of
%$J(P_{\gamma_*})$.
The idea is to approximate the boundary of $\Delta_{\gamma}$
with the iterates of the critical point $c_{\gamma}=-e^{2\pi i\gamma}/2$.
It is known that in this case the critical point itself is contained in the boundary.
The renormalization theory for golden-mean Siegel disks (constructed in \cite{McM}) implies
that the boundary $\Delta_{\gamma_*}$ is self-similar up to an 
exponentially small error. In particular, there exist  constants
$C>0$ and $\lambda>1$ such 
that
\begin{equation*}
d_H(\{P_{\gamma_*}^i(c_{\gamma_*}),i=0,\ldots, q_n\},\partial \Delta_{\gamma_*})<C\lambda^{-n}
\end{equation*}

Below we derive a similar estimate for all noble Siegel disks with 
constructive constants $C$ and $\lambda$.
For this, we do not need to invoke the whole power of renormalization theory. Rather,
we will use a theorem of Douady, Ghys, Herman, and Shishikura \cite{Do1} which specifically
applies to quadratic noble Siegel disks.

Noble (or more generally, bounded type) Sigel quadratic Julia sets may be 
constructed by means of quasiconformal surgery on a Blaschke product
$$f_\gamma(z)=e^{2\pi i\tau(\gamma)}z^2\frac{z-3}{1-3z}.$$
This map homeomorphically maps the unit circle $\TT$ onto itself with
a single (cubic) critical point at $1$. The angle
$\tau(\gamma)$ can be uniquely selected in such a way that the rotation number
of the restriction $\rho(f_\gamma|_\TT)=\gamma$.

For each $n$, the points $$\{1, f_\gamma(1),f^2_\gamma(1),\ldots,f^{q_{n+1}-1}_\gamma(1)\}$$
form the {\it $n$-th dynamical partition} of the unit circle. We have (cf. Theorem 3.1 of \cite{dFdM})
the following:

\begin{thm}[{\bf Universal real {\it a priori} bound}]
\label{real bound}
There exists an explicit constant $B>1$ independent of $\gamma$ and $n$
such that the following holds. Any two adjacent intervals $I$ and $J$
of the $n$-th dynamical partition of $f_\gamma$ are $B$-commensurable:
$$B^{-1}|I|\leq |J|\leq B|I|.$$
\end{thm}

\noindent
Let us now consider the mapping 
$\Psi$ which identifies the critical orbits of $f_\gamma$ and $P_\gamma$ by
$$\Psi:f^i_\gamma(1)\mapsto P^i_\gamma(c_\gamma).$$

\noindent
We have the following (Theorem 3.10 of \cite{YZ}):

\begin{thm}[{\bf Douady, Ghys, Herman, Shishikura}]
\label{surgery}
The mapping $\Psi$ extends to a $K$-quasiconformal homeomorphism of the plane
$\CC$ which maps the unit disk $\DD$ onto the Siegel disk $\Delta_\gamma$.
The constant $K$ depends on $B$ and $a_0,\ldots,a_k$  in a constructive fashion.
\end{thm}

Elementary combinatorics implies that each interval of the $n$-th dynamical
partition contains at least two intervals of the $(n+2)$-nd dynamical partition.
This in conjunction with \thmref{real bound} implies that the size of an interval
of the $(n+2)$-nd dynamical partition of $f_\gamma$ is at most $\tau^n$ where
$$\tau=\sqrt\frac{B}{B+1}.$$

We now  complete the proof of \lemref{bounded-type}. Denote $W_n$ the connected component
containing $0$ of the domain obtained by removing from the plane a closed disk of
radius $2K\tau^n$ around each point of 
$$\Omega_n=\{P_\gamma^i(c_\gamma),\;i=0,\ldots,q_{n+2}\}.$$
By Theorem \ref{surgery}, 
$$ \dist_H ( \Omega_n, \partial \Delta_\gamma) < K \tau^n,  $$
and we have $$W_n\subset \Delta_\gamma\text{ and }\dist_H(\partial \Delta_\gamma,\partial W_n)\leq 
\eps_n=2K\tau^n.$$
Any constructive algorithm for producing the Riemann mapping of a planar region 
(e.g. that of \cite{BB}) can be used to estimate the conformal radius $r(W_n,0)$ with
precision $\eps_n$. Denote this estimate $r_n$.

Elementary estimates imply that the Julia set $J(P_\gamma)\subset \overline{B(0,2)}$. By
Schwarz Lemma this implies 
$r(\Delta_\gamma,0)<2$.
By \lemref{rad-modulus1} we have 
$$|r(\Delta_\gamma,0)-r_n|
\le |r(\Delta_\gamma,0)- r(W_n,0)| + | r(W_n,0) - r_n | 
<4\sqrt{\eps_n}+\eps_n\underset{n\to\infty}{\longrightarrow}0,$$
and the proof is complete.

\section{Making Small Changes to $\Phi$ and $r$}
\label{sec changes}

%The next two 
%lemmas  are proved in \cite{BY}. We only reproduce their statements.

For a number $\gamma=[a_1,a_2,\ldots]\in\RR\setminus\QQ$ we denote
$$
\displaystyle\al_i(\gamma) = \frac{1}{a_i + \displaystyle\frac{1}{a_{i+1} + \displaystyle\frac{1}{a_{i+2}+\ldots}}}, 
$$
so that
$$
\Phi(\gamma) = \sum_{n\ge 1} \al_1(\gamma) \al_2(\gamma) \ldots \al_{n-1}(\gamma) \log \frac{1}{\al_n(\gamma)}.
$$

\noindent
We will show the following two lemmas.

\begin{lem}
\label{smlchg}
For any initial segment $I = (a_0, a_1, \ldots, a_n)$, write 
$\omega = [a_0, a_1, \ldots, a_n, 1, 1, 1, \dots]$. Then 
for any $\ve>0$, there is an $m>0$ and an integer $N$
such that if we write $\be = [a_0, a_1, \ldots, a_n, 1, 1, \ldots, 
1, N, 1, 1, \ldots]$, where the $N$ is located in the $n+m$-th 
position, then
$$
\Phi(\omega) + \ve < \Phi(\be) < \Phi(\omega) + 2 \ve. 
$$
\end{lem}

\begin{lem}
\label{notdeclem}
For  $\omega$ as above, for any $\ve>0$ there is an $m_0>0$, which 
can be computed from $(a_0, a_1, \ldots, a_n)$ and $\ve$,  
such that for any $m \ge m_0$, and for any tail 
$I = [a_{n+m}, a_{n+m+1}, \ldots]$ if we denote 
$$ \be^{I} = [a_1, a_2, \ldots, a_n, 1, 1, \ldots, 1, a_{n+m}, a_{n+m+1}, 
\ldots],$$ then
$$
\Phi(\be^{I}) > \Phi (\om) - \ve. 
$$
\end{lem}

We first prove lemma  \ref{smlchg}. Denote
$$
\Phi^{-} (\omega) = \Phi(\omega) - \al_0 (\om) \al_1 (\om) \ldots \al_{n+m-1} (\om) \log 
\frac{1}{\al_{m+n}(\om)}. 
$$
The value of the integer $m>0$ is yet to be determined.
Denote  $$\be^N = (a_0, a_1, \ldots, a_n, 1, 1, \ldots,
1, N, 1, 1, \ldots).$$ 

We will need the following estimates, which are proven by induction

\begin{lem}
\label{4lems}
For any $N$, the following holds:
\begin{enumerate}
\item 
\label{4l:p1}
For $i\le n+m$ we have $$\left| 
\log{\frac{\al_{i}(\be^{N})}{\al_{i}(\be^{N+1})}} \right|
< 2^{i-(n+m)}/N;$$
\item 
\label{4l:p2}
for $i < n+m$,
$$\left| 
\log{\frac{\al_{i}(\be^{N})}{\al_{i}(\be^{1})}} \right|
< 2^{i-(n+m)};$$
\item 
\label{4l:p3}
for $i< n+m$,$$
\left|
\log{\frac{\log \frac{1}{\al_{i}(\be^{N})}}{\log 
\frac{1}{\al_{i}(\be^{N+1})}}} \right|
< 2^{i-(n+m)+1};
$$
\item 
\label{4l:p4}
for $i < n+m-1$, $$
\left|
\log{\frac{\log \frac{1}{\al_{i}(\be^{N})}}{\log 
\frac{1}{\al_{i}(\be^{1})}}} \right|
< 2^{i-(n+m)+1}.
$$
\end{enumerate}
\end{lem}

The estimates yield the following. 

\begin{lem}
\label{lem5}
For any $\om$ of the form as in lemma \ref{smlchg} and for any $\ve >0$, 
there is an $m_0 >0$ such that for any $N$ and any $m \ge m_0$,
$$
 | \Phi^{-} (\be^{N}) -
\Phi^{-} ( \be^1) | < \frac{\ve}{4}.
$$
\end{lem}

\begin{proof}
The $\sum$ in the expression for $\Phi (\be^1)$ converges, 
hence there is an $m_1 >1$ such that the tail of the sum 
$\sum_{i\ge n+m_1} \al_1 \al_2 \ldots \al_{i-1} \log \frac{1}{\al_i} < 
\frac{\ve}{16}$. We will show how to choose $m_0 \ge m_1$ to
satisfy the conclusion of the lemma. 

We bound the influence of the change from $\be^1$ to 
$\be^N$ using lemma \ref{4lems} parts \ref{4l:p2} and \ref{4l:p4}. The influence on 
each of the ``head elements" ($i< n+m_1<n+m-1$) is bounded by 
$$
\left| \log{\frac{\al_1 (\be^1) \ldots \al_{i-1}(\be^1) \log 
\frac{1}{\al_i(\be^1)}}{\al_1 (\be^N) \ldots \al_{i-1}(\be^N) \log
\frac{1}{\al_i(\be^N)}}} 
\right| < \sum_{j=1}^{i-1} 2^{j-(n+m)} + 2^{i-(n+m)+1} < 2^{i -(n+m)+2}<
2^{m_1+2-m}.
$$
By making $m$ sufficiently large (i.e. by choosing a sufficiently large
$m_0$ we can ensure that 
$$
1- \frac{\ve}{16 \Phi(\be^1)} < \frac{\al_1 (\be^N) \ldots 
\al_{i-1}(\be^N) \log
\frac{1}{\al_i(\be^N)}}{\al_1 (\be^1) \ldots \al_{i-1}(\be^1) \log
\frac{1}{\al_i(\be^1)}} < 1 + \frac{\ve}{16 \Phi(\be^1)},
$$
hence
$$
\left| \al_1 (\be^N) \ldots
\al_{i-1}(\be^N) \log
\frac{1}{\al_i(\be^N)} - \al_1 (\be^1) \ldots
\al_{i-1}(\be^1) \log
\frac{1}{\al_i(\be^1)} \right| < 
$$
$$\frac{\ve}{16 \Phi(\be^1)}\al_1 (\be^1)
\ldots
\al_{i-1}(\be^1) \log
\frac{1}{\al_i(\be^1)}.
$$

Adding the inequality for $i=1,2, \ldots, n+m_1-1$ we obtain 
$$
\left| \sum_{i=1}^{n+m_1-1}  \al_1 (\be^N) \ldots
\al_{i-1}(\be^N) \log
\frac{1}{\al_i(\be^N)} - \sum_{i=1}^{n+m_1-1}  \al_1 (\be^1) \ldots
\al_{i-1}(\be^1) \log
\frac{1}{\al_i(\be^1)} \right| < 
$$
$$
\frac{\ve}{16 \Phi(\be^1)} \sum_{i=1}^{n+m_1-1}  \al_1 (\be^1) \ldots
\al_{i-1}(\be^1) \log
\frac{1}{\al_i(\be^1)} < \frac{\ve}{16 \Phi(\be^1)} \Phi(\be^1) = 
\frac{\ve}{16}.
$$
Hence the influence on the ``head" of $\Phi^-$ is bounded 
by $\frac{\ve}{16}$. 

To bound the influence on the ``tail" we consider three kinds of 
terms \\$\al_1 (\be^N) \ldots
\al_{i-1}(\be^N) \log
\frac{1}{\al_i(\be^N)}$: $n+m_1 \leq i \leq n+m-2$, $i=m+n-1$ and
$i \ge m+n+1$ (recall that $i=n+m$ is not in $\Phi^-$).

{\bf For $n+m_1 \leq i \leq n+m-2$:} 
$$
\left| \log{\frac{\al_1 (\be^1) \ldots \al_{i-1}(\be^1) \log
\frac{1}{\al_i(\be^1)}}{\al_1 (\be^N) \ldots \al_{i-1}(\be^N) \log
\frac{1}{\al_i(\be^N)}}}
\right| < \sum_{j=1}^{i-1} 2^{j-(n+m)} + 2^{i-(n+m)+1} < 2^{i -(n+m)+2}\le
1.
$$
Hence in this case each term can increase at most by a factor of $e$. 

{\bf For $i=n+m-1$} Note that the change decreases 
$\log{\frac{1}{\al_{n+m-1}}}$
so that $\log{\frac{1}{\al_{n+m-1}(\be^N)}} \le 
\log{\frac{1}{\al_{n+m-1}(\be^1)}}$, hence we have 
$$
\log{\frac{\al_1 (\be^N) \ldots \al_{i-1}(\be^N) \log
\frac{1}{\al_i(\be^N)}}{\al_1 (\be^1) \ldots \al_{i-1}(\be^1) \log
\frac{1}{\al_i(\be^1)}}} \le 
\log{\frac{\al_1 (\be^N) \ldots \al_{i-1}(\be^N)}{\al_1 (\be^1) \ldots 
\al_{i-1}(\be^1)}} < 
$$
$$
\sum_{j=1}^{n+m-2} 2^{j-(n+m)} < \frac{1}{2}. 
$$
Hence this term could increase by a factor of $\sqrt{e}$ at most. 

{\bf For $ i \geq n+m+1$:}
Note that $\al_j$ for $j>n+m$ are not affected by the change, 
and the change decreases $\al_{n+m}$, so that $\al_{n+m} (\be^N) \le 
\al_{n+m} (\be^1)$. 
Hence 
$$
\log{\frac{\al_1 (\be^N) \ldots \al_{i-1}(\be^N) \log
\frac{1}{\al_i(\be^N)}}{\al_1 (\be^1) \ldots \al_{i-1}(\be^1) \log
\frac{1}{\al_i(\be^1)}}} = 
\log{\frac{\al_1 (\be^N) \ldots \al_{n+m}(\be^N)}{\al_1 (\be^1) \ldots 
\al_{n+m}(\be^1)}} \le
$$
$$
\log{\frac{\al_1 (\be^N) \ldots \al_{n+m-1}(\be^N)}{\al_1 (\be^1) \ldots
\al_{n+m-1}(\be^1)}} < \sum_{j=1}^{n+m-1} 2^{j-(n+m)} < 1
$$
So in this case each term could increase by a factor of $e$ 
at most. 

We see that after the change each term of the tail could increase 
by a factor of $e$ at most. The value of the tail remains 
positive in the interval $(0, \frac{e \ve}{16}]$, hence the 
change in the tail is bounded  by $\frac{e \ve}{16} < \frac{ 3 \ve}{16}$. 

So the total change in $\Phi^-$ is bounded by
$$
\mbox{change in the ``head"}~~ +~~
\mbox{change in the ``tail"} ~~ < \frac{\ve}{16} + \frac{3 \ve}{16} = 
\frac{\ve}{4}.
$$
\end{proof}

Lemma \ref{lem5} immediately yields:

\begin{lem}
\label{lem6}
For any $\ve$ and for the same $m_0 (\ve)$ as in lemma \ref{lem5}, 
for any $m \ge m_0$ and $N$, 
$$
 | \Phi^{-} (\be^{N}) -
\Phi^{-} ( \be^{N+1}) | < \frac{\ve}{2}.
$$ 
\end{lem}

Denote $\Phi^{1} (\omega) =  \al_0 (\om) \al_1 (\om) \ldots \al_{n+m-1} (\om) \log 
\frac{1}{\al_{m+n}(\om)} = \Phi(\om) - \Phi^{-} (\om)$. We are now ready to prove the following. 

\begin{lem}
\label{lem9}
For sufficiently large $m$, for any $N$, 
$$
\Phi^1 (\be^{N+1}) - \Phi^1 (\be^N) < \frac{\ve}{2}.
$$
\end{lem}

\begin{proof}
According to lemma \ref{4lems} part \ref{4l:p1} we have 
$$
\left|
\log{\frac{\al_{1}(\be^{N+1})\ldots 
\al_{n+m-1}(\be^{N+1})}{\al_{1}(\be^{N}) \ldots \al_{n+m-1}(\be^{N})}} 
\right| < \sum_{i=1}^{n+m-1} 2^{i-(n+m)}/N < \frac{1}{N}. 
$$
Hence $\al_{1}(\be^{N+1})\ldots
\al_{n+m-1}(\be^{N+1}) < \al_{1}(\be^{N}) \ldots \al_{n+m-1}(\be^{N}) 
e^{1/N}$, and 
$$
\Phi^1 (\be^{N+1}) < \Phi^1 (\be^N) e^{1/N}
\frac{\log{\frac{1}{\al_{n+m}(\be^{N+1})}}}
     {\log{\frac{1}{\al_{n+m}(\be^{N})}}} = 
\Phi^1 (\be^N) e^{1/N}
\frac{\log (N+1+\phi)}{\log(N+\phi)}.
$$
Hence 
$$
\Phi^1 (\be^{N+1}) - \Phi^1 (\be^N) < \Phi^1 (\be^N) \left(e^{1/N}
\frac{\log (N+1+1/\phi)}{\log(N+1/\phi)}-1\right) <
$$
$$
\Phi^1 (\be^N)\left((1+\frac{e}{N}) \frac{\log 
(N+1+1/\phi)}{\log(N+1/\phi)}-1\right).
$$
We make the following calcualtions.
Denote $x = \frac{\log
(N+1+1/\phi)}{\log(N+1/\phi)}$, then $(N+1/\phi)^x = N+1+1/\phi$,
$(N+1/\phi)^{x-1} = 1 + \frac{1}{N+1/\phi}< e^{\frac{1}{N+1/\phi}}$. 
$N+1/\phi > e^{1/3}$, and so $x-1 < \frac{3}{N+1/\phi} < \frac{3}{N}$,
thus $x < 1 + \frac{3}{N}$. 

It is not hard to see that $\al_{k-1} \al_k < 1/2$ for all $k>1$, and we have 
$$
\Phi^1 (\be^N) = \al_1 (\be^N) \ldots \al_{n+m-1} (\be^N) \log 
\frac{1}{\al_{n+m} (\be^N)} < \left( \frac{1}{2} \right)^{(n+m-1)/2} 
\log ( N + 1/\phi). 
$$
Thus 
$$
\Phi^1 (\be^{N+1}) - \Phi^1 (\be^N) <
\Phi^1 (\be^N)\left((1+\frac{e}{N}) \frac{\log
(N+1+1/\phi)}{\log(N+1/\phi)}-1\right) <
$$
$$
\left( \frac{1}{2} \right)^{(n+m-1)/2}
\log ( N + 1/\phi) \left( (1+e/N)(1+3/N) -1 \right) < 
\left( \frac{1}{2} \right)^{(n+m-1)/2} \log ( N + 1/\phi)
\frac{14}{N}.
$$
Since $\frac{14}{N} \in o(1/\log ( N + 1/\phi))$, this expression 
can be always made less than $\frac{\ve}{2}$ by choosing $m$ 
large enough. 
\end{proof}

Since $\Phi = \Phi^- + \Phi^1$, summing the inequalities  in Lemmas 
\ref{lem6} and \ref{lem9} yields the following.

\begin{lem}
\label{lem10}
For sufficiently large $m$, for any $N$, 
$$
\Phi (\be^{N+1}) - \Phi (\be^N) < \ve.
$$
\end{lem}

It is immediate from the formula of $\Phi(\beta^N)$ that:

\begin{lem}
\label{lem11}
$$
\lim_{N \rightarrow \infty} \Phi (\be^N) = \infty.$$
\end{lem}

We are now ready to prove lemma \ref{smlchg}. 

\begin{proof} (of lemma \ref{smlchg}). Choose $m$ large enough for lemma 
\ref{lem10} to hold. Increase $N$ by one at a time starting with $N=1$. 
We know that $\Phi(\be^1) = \Phi(\om) < \Phi(\om) + \ve$, and 
by lemma \ref{lem11}, there exists an $M$ with $\Phi(\be^M) > \Phi(\om) + 
\ve$. Let $N$ be the smallest such $M$.  Then $\Phi(\be^{N-1}) \le
\Phi(\om) +\ve$, and by lemma \ref{lem10}
$$
\Phi(\be^N) < \Phi(\be^{N-1}) + \ve \le \Phi(\om) + 2 \ve.
$$
Hence 
$$
\Phi(\om) +\ve < \Phi(\be^N) <  \Phi(\om) + 2 \ve.
$$
Choosing $\be = \be^N$ completes the proof. 
\end{proof}

\noindent
The second part of the following Lemma follows by the same argument as Lemma \ref{lem5} 
by taking $N\ge 1$ to be an arbitrary real number, not necessairily
an integer. The first part is obvious, since the tail of $\om$ has only $1$'s. 

\begin{lem}
\label{auxbdtl1}
For an $\om = \be^1$ as above, for any $\ve>0$ there is an $m_0 >0$, 
such that for any $m \ge m_0$, and for any tail 
$I = [a_{n+m}, a_{n+m+1}, \ldots]$ if we denote 
$$ \be^{I} = [a_1, a_2, \ldots, a_n, 1, 1, \ldots, 1, a_{n+m}, a_{n+m+1}, 
\ldots],$$ then
$$
\sum_{i \ge n + m} \al_1 (\be^1) \al_2 (\be^1) \ldots  \al_{i-1} (\be^1)
\log \displaystyle\frac{1}{\al_i(\be^1)} < \ve, 
$$
and 
$$
\sum_{i=1}^{n+m-1} \left| \al_1 (\be^I) \ldots
\al_{i-1}(\be^I) \log
\displaystyle\frac{1}{\al_i(\be^I)} - \al_1 (\be^1) \ldots
\al_{i-1}(\be^1) \log
\displaystyle\frac{1}{\al_i(\be^1)} \right| < 
\ve.$$
\end{lem}

We can now prove lemma \ref{notdeclem}.

\begin{proof} (of lemma \ref{notdeclem}). 
Applying lemma \ref{auxbdtl1} with $\frac{\ve}{2}$ instead of $\ve$, we get
$$
\Phi(\be^I) - \Phi(\om) = \sum \{ \mbox{``head"($\be^I$)}-\mbox{``head"($\om$)}\} + 
\sum \{ \mbox{``tail"($\be^I$)}-\mbox{``tail"($\om$)}\} > 
$$
$$
-\frac{\ve}{2}-\sum \{\mbox{``tail"($\om$)}\} > -\frac{\ve}{2}-\frac{\ve}{2} = -\ve. 
$$
\end{proof}

\comm{
%%%%\subsection{Making the tail $[1, 1, 1, \ldots ]$.} 
\noindent
We will also need the 
following Lemma in the proof of the Main Theorem.

\begin{lem}
\label{tailto1}
Let $\al = [a_0, a_1, a_2, \ldots]$ and let $\ve>0$ be given. 
Then there is an $N = N(\ve)$ such that for any $n \ge N$ 
we have $\Phi(\al_n)<\Phi(\al)+\ve$, where 
$\al_n = [a_0, a_1, \ldots, a_{n-1}, 1, 1, 1, \ldots]$. 
\end{lem}
}

We will need a {\em computable} version of Lemma \ref{smlchg}
for modifying the conformal radius of the corresponding 
Julia set. 

\begin{lem}
\label{smlchgr}
For any given initial segment $I = (a_0, a_1, \ldots, a_n)$ and $m_0>0$, write 
$\omega = [a_0, a_1, \ldots, a_n, 1, 1, 1, \dots]$. Then 
for any $\ve>0$, we can uniformly compute  $m>m_0$ and an integer $N$
such that if we write $\be = [a_0, a_1, \ldots, a_n, 1, 1, \ldots, 
1, N, 1, 1, \ldots]$, where the $N$ is located in the $n+m$-th 
position, we have 
\beq
\label{rcond}
r(\omega) - 2 \ve < r(\be) < r(\omega) -  \ve,
\eeq
and
\beq
\label{phiinc}
\Phi(\be)> \Phi(\om).
\eeq
\end{lem}

\begin{proof}
We first show that such $m$ and $N$ {\em exist}, and then 
give an algorithm to compute them. By Lemma \ref{smlchg}
we can increase $\Phi(\omega)$ by any controlled amound 
by modifying one term arbitrarily far in the expansion. 

By Theorem \ref{phi-cont}, $f: \theta\mapsto \Phi(\theta)+\log r(\theta)$
extends to a continuous function. Hence for any $\ve_0$ there is a $\de$ such 
that $|f(x)-f(y)|<\ve_0$ whenever $|x-y|<\de$. In particular, 
there is an $m_1$ such that $|f(\be)-f(\om)| < \ve_0$ whenever $m \ge m_1$.

This means that if we choose $m$ large enough, a controlled increase 
of $\Phi$ closely corresponds to a controlled drop of $r$ by a corresponding 
amount, hence there are $m>m_0$ and $N$ such that 
\eref{rcond} holds. \eref{phiinc} is satisfied almost automatically.
The only problem is to {\em computably} find such $m$ and $N$. 

To this end, we apply Lemma \ref{bounded-type}. It implies that for any specific $m$ 
and $N$ we can compute $r(\be)$. This means that we can find the suitable 
$m$ and $N$, by enumerating all the pairs $(m,N)$ and exhaustively checking 
\eref{rcond} and \eref{phiinc} for all of them. We know that eventually we will find a 
pair for which \eref{rcond} and \eref{phiinc} hold. 
\end{proof}

\section{Proving the Main Theorem}
\label{sec:proof}

There are countably many oracle Turing Machines. Let us enumerate  them in some arbitrary
computable fashion $M_1^{\phi}, M_2^{\phi}, \ldots$ so that every machine appears 
infinitely many times in the enumeration. 
Recall that $r(\theta)$ is the conformal radius of the Siegel disk associated
with the polynomial $P_\theta (z) = z^2 + e^{2 \pi i \theta} z$, or zero,
if $\theta$ is not a Brjuno number.

We will argue by induction.
On each iteration $i$ of the argument we shall 
maintain an initial segment $I_i = [a_0, a_1, \ldots, a_{N_i}]$ 
an  interval $H_i=[l_i, r_i]$, and $\ell_i=\ell(H_i)=r_i-l_i$  such that the following properties are maintained: 
\beq
\label{cond1}
r_i = r(\ga_i), \mbox{   where  } \ga_i = [I_i, 1, 1, \ldots],
\eeq
and
\beq 
\label{cond2}
\mbox{For any } \be = [I_i, t_{N_i+1}, t_{N_i+2}, \ldots] \mbox{ with } 
r(\be) \in 
[l_i, r_i], 
\eeq

\noindent
the machine  $M_i^{\phi}$ requires at least the time $h(2 \left\lceil -\log \ell_i\right\rceil+1)$ to 
compute the $\frac{\ell_i^2}{2}$-approximation to $J(P_\be)$. And

\beq
\label{cond3}
\mbox{for $i\geq 1$, 
$\Phi(\beta)>\Phi(\gamma_{i-1})-2^{-(i-1)}$, for any  $\be = [I_i, 
t_{N_i+1}, t_{N_i+2}, \ldots]$.}
\eeq

\noindent
Moreover, the intervals we construct are nested: $[l_i,r_i]\subset[l_{i-1},r_{i-1}]$, and the sequence
$I_i$ contains $I_{i-1}$ as the initial segment. 
The numbers $2 \left\lceil -\log \ell_i\right\rceil$ form a strictly increasing sequence.

\comm{
\noindent
We also keep track of  the {\em  variation} 
\beq 
\label{ivar}
T_j =  \sum_{i \ge 0} \left| 
\al_1 (\ga_{i+1}) \al_2 (\ga_{i+1}) \ldots \al_{j-1} (\ga_{i+1}) \log \frac{1}{\al_j (\ga_{i+1})} -
\al_1 (\ga_{i}) \al_2 (\ga_{i}) \ldots \al_{j-1} (\ga_{i}) \log \frac{1}{\al_j (\ga_{i})}\right|
\eeq
at each step of the process, in order to be able to pass to the limiting sequence $[a_0,a_1,a_2,a_3,\ldots]$ 
in the end.

\medskip}
\noindent
%We show how to initialize and maintain these properties by induction.
%
For the basis of induction, set $I_0 = [1]$, $r_0= r(\ga_0)<2$ (by the Schwarz Lemma) 
and $l_0 = r_0/2$, where
$\ga_0 = [1,1,1, \ldots]$. Then for $i=0$ condition \eref{cond1} holds by definition
and conditions \eref{cond2} and \eref{cond3} hold because they are empty.

\subsection*{The induction step.} We now have the conditions \eref{cond1}, 
\eref{cond2}
and \eref{cond3} for some $i$ and would like to extend them
to $i+1$. 

\medskip

Consider the machine $M_{i+1}^{\phi}$. Set $\ell_{i+1} = \frac{\ell_{i}}{20}$.
Simulate $M_{i+1}^{\phi}$ on $\ga_i$ for at most $h(2 \left\lceil -\log \ell_{i+1}\right\rceil+1)$
 steps to compute $J_{\ga_i}$ with precision $\frac{\ell_{i+1}^2}{2}$. The machine 
 reads at most $h(2 \left\lceil -\log \ell_{i+1}\right\rceil+1)$ bits of the input, and 
 we can compute  $m_0$ such that this run does not distinguish between $\ga_i$ and
 $\ga = [I_i, 1, 1, \ldots,1, N_{m_0+1}, N_{m_0+2}, \ldots]$.
 There are two cases:

 \medskip
\noindent
{\bf Case 1:} $M_{i+1}^{\phi}$ does not terminate in the assigned time, or
 does not output a proper set. In this case,
 we proceed by setting $I_{i+1} = [I_i, 1, \ldots, 1]$ (with $1$'s up to position $m_0$), $\ga_{i+1}=\ga_i$, 
 $r_{i+1} = r_i$, and $l_{i+1} = r_{i+1}-\ell_{i+1}$. By lemma \ref{notdeclem}, we can choose 
 sufficiently many $1$'s in $I_{i+1}$, so that for any $\be$ beginning with 
 $I_{i+1}$, we have $\Phi(\be)> \Phi(\ga_{i})-2^{-i}$. 
 
\medskip

\noindent
 {\bf Case 2:} $M_{i+1}^{\phi}$ outputs a set $S$. Compute the conformal radius $r(S)$. 
Considerations of  Schwarz Lemma imply that for any  quadratic Siegel disk,
$r(\Delta)<2$. Using the above consideration to bound the constant in 
Lemma \ref{rad-modulus1}, we know that for any Julia set $J({P_\om})$ which is $\ell_{i+1}^2$-accurately 
described by $S$, we have $|r(\om)-r(S)|<4\sqrt{2} \ell_{i+1}<6 \ell_{i+1}$. Again, there are 
two cases (if both hold, it doesn't matter which way to proceed):

\noindent
{\bf Subcase 2a:} $r_i - \ell_{i+1} > r(S)+ 8 \ell_{i+1}$. In this case
we proceed by setting $I_{i+1} = [I_i, 1, \ldots, 1]$ (with $1$'s up to position $m_0$), $\ga_{i+1}=\ga_i$, 
 $r_{i+1} = r_i$, and $l_{i+1} = r_{i+1}-\ell_{i+1}$. 

\noindent
{\bf Subcase 2b:} $l_i + 2 \ell_{i+1} < r(S)- 8 \ell_{i+1}$. By Lemma 
\ref{smlchgr}, we can select $\ga_{i+1}=
[I_{i+1},1,1, \ldots]$ by  modifying $\ga_i$ at an arbitrarily far 
position, and set $r_{i+1}=r(\ga_{i+1})$
so that $\Phi(\ga_{i+1})>\Phi(\ga_i)$, $r_{i+1}\ge l_i+\ell_{i+1}$ and  
$[r_{i+1}-\ell_{i+1},r_{i+1}] \cap [r(S)-8 \ell_{i+1}, r(S)
+8 \ell_{i+1}]=\emptyset$. The number $r_{i+1}$ is computable
since it is the conformal radius of a noble Siegel disk.
Set $l_{i+1}= r_{i+1}-\ell_{i+1}$. We see that the induction 
is maintained for these parameters. 

\smallskip 
\noindent
In either subcase, by Lemma \ref{notdeclem}, we can add 
 sufficiently many $1$'s to $I_{i+1}$, so that for any $\be$ beginning with 
 $I_{i+1}$, we have $\Phi(\be)> \Phi(\ga_{i})-2^{-i}$, and condition 
\eref{cond3} is satisfied.

\medskip

\comm{Choose an $\eps>0$ such that $(r_i-l_i)/(i+1)>\eps$. By the Koebe Theorem there exists 
$k\in\NN$ such that the following
holds:\\[6pt]
 Suppose $r(\theta_1)\in[l_i,r_i]$  and $|r(\theta_2)-r(\theta_1)|>\eps$. Then 
$$d_H(J(P_{\theta_1}),J(P_{\theta_2}))>2^{-(k-1)}.$$

\medskip
\noindent
Recall the definition of $\gamma_i$ in \eref{cond1}, and
``feed'' the first $h(k)$ digits of $\gamma_i$ to each of the machines 
$M_1^\phi,\ldots,M^\phi_{i+1}$. That is,
set the value of the oracle $\phi$ to coincide with $\gamma_i$ on the first $k$ binary digits, and denote
$K_j$ the output of the TM $M^\phi_j$ with this oracle {\it in time $h(k)$}. It is possible, of course,
that $K_j=\emptyset$. By the Pigeonhole Principle, there exists a subinterval $[a,b]\subset[l_i,r_i]$ 
of length at least $\eps$ with the
following property:\\[6pt]
 Suppose that $\be\in\TT$ has the same $h(k)$ first digits as $\gamma$,
 and $r(\beta)\in[a,b]$.
Then 
$$d_H(J(P_{\be}),K_j)>2^{-(k-1)}\text{ for all }1\leq j\leq i+1.$$

\medskip
\noindent
We now select $k$ so that in \lemref{auxbdtl1} for $\omega=\gamma_i$ and $m_0=k$ we have $\eps<2^{-i}$.
Set $A=[I_i,\underbrace{1,1,\ldots,1}_{s}]$, where $s$ is such that $A$ has the first $h(k)$ digits 
in common with $\gamma_i$.  Increasing $s$ if necessary,
by \lemref{smlchg}, for the parameter  $\be=[A,\underbrace{1,1,\ldots,1}_{s},N,1,1,\ldots]$ we have
 $r(\be)\in[a,b]$. 
Set $I_{i+1}=[A,1,\underbrace{1,1,\ldots,1}_{s},N]$, 
$\ell_{i+1}=k$, $[l_{i+1},r_{i+1}]=[a,r(\beta)]$.
}

%\noindent
%A bound on the variation $T_j$ is necessary to conclude:

\begin{lem}
Denote $\ga = \lim_{i \rightarrow \infty} \ga_i$. Then the following equalities hold:
$$
\Phi(\ga) = \lim_{i \rightarrow \infty}\Phi(\gamma_i)~~~~\mbox{$~$and$~$}~~~~
r(\ga) = \lim_{i \rightarrow \infty}r(\gamma_i).
$$
\end{lem}

\begin{proof}
By the construction, the limit $\ga = \lim \ga_i$ exists. We also 
know that the sequence $r(\ga_i)=r_i$ converges uniformly to 
some number $r$, and that the sequence $\Phi(\ga_i)$ is 
monotone non-decreasing, and hence converges to a value $\psi$ ({\em a priori}
we could have $\psi=\infty$). By the Carath{\'e}odory Kernel Theorem (see e.g. \cite{Pom}), we have $r(\ga) \ge r >0$, so 
$\psi < \infty$.
On the other hand, by the property we have maintined through 
the construction, we know that $\Phi(\ga)>\Phi(\ga_i)-2^{-i}$ for all $i$. 
Hence $\Phi(\ga) \ge \psi$. 

From \cite{BC} we know that 
\beq
\label{BCCont}
\psi + \log r= \lim (\Phi(\ga_i) + \log r(\ga_i)) = \Phi(\ga) + \log r(\ga).
\eeq
Hence we must have $\psi = \Phi (\ga)$, and $r = r(\ga)$, which completes
the proof. 
\end{proof}

The conformal radius 
$r(\ga)$ is computable, since the convergence $r(\ga_i)\rightarrow 
 r(\ga)$ is uniform. Thus $J_{P_\ga}$ is also computable by \thmref{thm BBY}. 
By construction, it satisfies all of the required properties.
Note that the value $\ga$ itself is also computable.

\comm{

There exist $\ve_0>0$ and $m_0\in\NN$ such that for every $\be\in S$ of the form 
$$\be=[I_i,\underbrace{1,1,\ldots,1}_{m_0},\ldots]\text{ we have }|r_i-r(\be)|>\ve_0.$$
In this case, select $0<\ve\leq\ve_0$ such that $r_i-\ve>l_i$.
Set
 $$I_{i+1}=[I_i,\underbrace{1,1,\ldots,1}_{m_0}],\;l_{i+1} = r_i - \ve,
r_{i+1} = r_i,\text{ and }\ell_{i+1}=\ell_i.$$ 
Since 
$\ga_{i+1} = [I_i,1, 1, \ldots ] = \ga_i$, we have
$r(\ga_{i+1}) = r(\ga_i) = r_i = r_{i+1}$ and the condition 
\eref{cond1} is satisfied. 

Suppose $\be = [I_{i+1}, t_{N_{i+1}+1}, t_{N_{i+1}+2}, \ldots]$ with $r(\be) \in 
[l_{i+1}, r_{i+1}]$. $I_{i+1}$ is an extension of $I_i$, and 
the machines $M_1^{\phi}, M_2^{\phi}, \ldots, M_i^{\phi}$ cannot compute $2^{-\ell_{i+1}}$-approximation
of $J(P_\beta)$ in time less than $h(\ell_{i+1})$ by the induction assumption.
on $\be$ by the induction assumption. $M_{i+1}^\phi$ fails 
on $\be$ because $r(\be) \in [l_{i+1}, r_{i+1}]$, and so 
$r(\be) \notin R$. This shows \eref{cond2} and completes
the proof for this case. 
Note that in this case we have $\ga_{i+1}=\ga_i$, and so $T_i=0$. 

\medskip
\noindent
The complementary case is the main part of the argument:

\noindent
{\bf Case 2.} 
For every $\ve>0$ and $m\in\NN$ we can find $\beta\in S$ starting with $I_i$ followed by $m$ ones
so that
\begin{equation}
\label{be1}
r_i - \ve< r(\be) \le r_i
\end{equation}
Choose the number  $m$ of $1$'s to be such that  \lemref{auxbdtl1}
holds with $\ve=2^{-i}$.

Choose an $\ve$ such that $r_i - 3 \ve > l_i$, and let $\be \in S$ as above 
satisfying $r_i - \ve< r(\be) \le r_i$. Denote 
$$\ve_0 = \min \left( \displaystyle\frac{\log(r_i-\ve)-\log(r_i- 2 \ve)}{8},
\frac{\log(r_i-2 \ve)-\log(r_i- 3 \ve)}{8}, 2^{-i} \right)>0.$$

\noindent
Theorem \ref{phi-cont} says that $\theta \mapsto \Phi(\theta)+\log r(\theta)$ 
extends to a continuous $1$-periodic function on $\RR$. Denote this function by
$f(\theta)$. $f$ is continuous and periodic, hence it must be uniformly
continuous, and there is a $\de>0$ such that if $|x-y|< \de$ then $|f(x)-f(y)|<\ve_0$. 

\noindent
By \thmref{cont2} the Julia sets of $P_\theta$ depend continuously on $\theta$ with respect to the
Hausdorff metric for $\theta\in S$.
By Proposition \ref{radius-continuous} the conformal 
radius $r(\bullet)$ is continuous on $S$. Hence there is a
$\de_0>0$ such that $$|  r(x) -  r(\be)| < \ve\text{ whenever } 
|x - \be| < \de_0\text{ and }x \in S.$$

We choose $m$ large enough, so that for any $\zeta = [I_i, \underbrace{1,1, \ldots, 1}_{m}, \ldots]$, 
$| \ga_i - \zeta |< \de$.
Write $\be = [I_i, \underbrace{1,1, \ldots, 1}_{m},  t_{N_i+1}, t_{N_i+2}, \ldots]$. By Lemma \ref{tailto1}, there is 
an $N$ such that for any $n \ge N$, $\be_n = [I_i, \underbrace{1,1, \ldots, 1}_{m}, t_{N_i+1}, \ldots, t_{N_i+n}, 1, 1, \ldots]$
satisfies $\Phi(\be_n) < \Phi(\be) + \ve_0$.  We can choose $n \ge N$ large 
enough so that for any $x$ with the initial segment 
$I_i' = [I_i, \underbrace{1,1, \ldots, 1}_{m}, t_{N_i+1}, t_{N_i+2}, \ldots , t_{N_i+n}]$, $|x - \be| < \de/2$ 
and $|x - \be| < \de_0$.

\smallskip

Start with $\om_0 = \be_n = [I_i', 1, 1, \ldots]$. We have 
$| \om_0 - \be | < \de$, hence $|f(\om_0)-f(\be)|< \ve_0$. So

$$
\log r(\om_0) = f(\om_0) - \Phi(\om_0 ) > f(\be) - \ve_0 - \Phi(\be)- \ve_0= 
\log r(\be) - 2 \ve_0.
$$

\noindent
 By Lemma \ref{smlchg} we can 
extend $\om_0$ to $\om_1= [I_i^1, 1, 1, \ldots]$ so that $\Phi(\om_0)+2 \ve_0< 
\Phi(\om_1) < \Phi(\om_0)+4 \ve_0$. Hence 
$$
\log (r(\om_1)) = f(\om_1) - \Phi(\om_1) > f(\om_0) - \ve_0 - \Phi(\om_0)
- 4 \ve_0 = \log (r(\om_0)) - 5 \ve_0,
$$
and
$$
\log (r(\om_1)) = f(\om_1) - \Phi(\om_1) < f(\om_0) + \ve_0 - \Phi(\om_0)
- 2 \ve_0 = \log (r(\om_0)) -  \ve_0.
$$
Hence $\log (r(\om_0)) -   5 \ve_0 < \log (r(\om_1)) < \log (r(\om_0)) -  \ve_0$. 
In the same fasion we can extend $I_i^1$ to $I_i^2$, and obtain 
$\om_2= [I_i^2, 1, 1, \ldots]$ so that 
$\log (r(\om_1)) -   5 \ve_0 < \log (r(\om_2)) < \log (r(\om_1)) -  \ve_0$.

\medskip
Recall that $\log (r(\om_0)) > \log (r(\be))- 2 \ve_0 > \log (r_i -  \ve) - 2 \ve_0
\ge  \log (r_i -  2 \ve) + 6 \ve_0$. Hence
after finitely many steps we will obtain $I_i^k$ and $\om_k= [I_i^k, 1, 1, \ldots]$
such that 
$$
\log ( r_i - 3 \ve ) + \ve_0 < \log ( r(\om_k)) < \log (r_i - 3 \ve ) + 6 \ve_0 < 
\log(r_i - 2 \ve). 
$$
Choose $I_{i+1} = I_i^k$, $\ga_{i+1} = \om_k$,  $l_{i+1}=l_i$, and $r_{i+1} = r(\om_k)$. We have 
$l_{i+1} < r_{i+1} < r_i - 2 \ve$. Condition \eref{cond1} is satisfied by definition.
Condition \eref{cond2} is satisfied for the first $i$ machines because $[l_{i+1}, r_{i+1}]$ 
is a subinterval of $[l_i, r_i]$, and $I_{i+1}$ is an extension of $I_i$. Condition 
\eref{cond2} is satisfied fot the $i+1$-st machine $M$, because $M$  could work correctly
on $x = [I_i', \ldots]$ only with  $r(x) \in [r_i-2 \ve, r_i + \ve]$, 
and $[l_{i+1}, r_{i+1}]$ is disjoint from  $[r_i-2 \ve, r_i + \ve]$.

\medskip
We now estimate the total variation $T_i$ in the 
process described above. By Lemma \ref{auxbdtl1}
we know that 
$$
\sum_{j \ge N_i+1} \al_1 (\ga_i) \al_2 (\ga_i) \ldots  \al_{j-1} (\ga_i)
\log \frac{1}{\al_j (\ga_i)} < 2^{-i}, 
$$
and 
$$
\sum_{j=1}^{N_i} \left| \al_1 (\ga_{i+1}) \ldots
\al_{j-1}(\ga_{i+1}) \log
\frac{1}{\al_j(\ga_{i+1})} - \al_1 (\ga_i) \ldots
\al_{j-1}(\ga_i) \log
\frac{1}{\al_j(\ga_i)} \right| < 
2^{-i}.$$
Denote $d = \Phi(\ga_{i+1}) - \Phi(\ga_i)$. 
%We know that $d \ge 0$ by the construction. 
Refer to  the first $N_i$ terms as the ``head terms", 
and the rest of the term as the ``tail terms". Denote the total increase 
in head terms from $\ga_i$ to $\ga_{i+1}$ by $h^{+}$, and the total 
decrease by $h^{-}$. Similarly define $t^{+}$ and $t^{-}$. Then
$d = h^+ - h^- + t^+ - t^-$, and $T_i = h^+ + h^- + t^+ + t^-$.
By the second inequality we have $h^+ + h^- < 2^{-i}$. 
By the first one we have $t^- < 2^{-i}$. Hence 
$$t^+ \le d + h^+ + h^- + t^- < d + 2 \cdot 2^{-i},$$ and 
$T_i < 2^{-i} + 2^{-i} + d + 2 \cdot 2^{-i} = d + 4 \cdot 2^{-i}$. We now 
bound $d$. 
$$
d = \Phi(\ga_{i+1}) - \Phi(\ga_i) = f(\ga_{i+1}) - \log ( r(\ga_{i+1})) -
f(\ga_{i}) + \log ( r(\ga_{i})) < $$ $$ \ve_0 +  \log ( r(\ga_{i})) -  \log ( r(\ga_{i+1}))<
2^{-i} + \log (r_i) - \log (r_{i+1}). 
$$
Hence $T_i < \log (r_i) - \log (r_{i+1}) +  5 \cdot 2^{-i}$. 

Overall we get a bound on the sum of {\em all} the $T_i$'s:
$$
\sum_{i=0}^{\infty} T_i < \log (r_0) - \log (l_0) + 5 \cdot \sum_{i=1}^{\infty} 2^{-i} = 
 \log (r_0) - \log (l_0) + 10 < \infty.
$$

\subsection*{Passing to the Limit}
The completion of the proof relies on the following Lemma:

\begin{lem}
\label{main2}
Let $I_i$ be the initial segment of $\ga_i$ as above. Then the sequence 
$\{ \Phi(\ga_i) \}$ converges to a limit 
$$
L = \lim_{i \rightarrow \infty} \Phi(\ga_i).
$$ 
Denote by $I$ the inductive limit of the segments above, 
and let $\ga$ be the number represented by $I$, then 
$\Phi(\ga) = L$.
\end{lem}

\begin{proof}
Denote $b_{ij} = \al_1 (\ga_j) \al_2 (\ga_j) \ldots \al_{i-1} (\ga_j) \log \frac{1}{\al_{i} (\ga_j)}$, 
and  $b_i =  \al_1 (\ga) \al_2 (\ga) \ldots \al_{i-1} (\ga) \log \frac{1}{\al_{i} (\ga)}$.
Our first claim is that for any $i$, $\lim_{j \rightarrow \infty} b_{ij} = b_i$. 

Note that on step $j$ the position being modified is at least $j$, hence
by Lemmas \ref{lem2} and \ref{lem4} we see that for any $j>i$, 
$$
\left| \log \frac{\al_i ( \ga_j)}{\al_i (\ga)} \right| < 2^{i-j}, \mbox{   and   }
\left| \log \frac{\log \frac{1}{\al_i ( \ga_j)}}{\log \frac{1}{\al_i (\ga)}} \right| < 2^{i-j+1}.
$$
Hence for all $j>i$,
$$
\left| \log \frac{b_{ij}}{b_i} \right| < \sum_{k=1}^{i-1} 2^{k-j} + 2^{i-j+1} < 
2^{i-j} + 2^{i-j+1} < 2^{i-j+2}. 
$$ 
Hence $\lim_{j \rightarrow \infty} \frac{b_{ij}}{b_i} = 1$. We know that $b_i$ is finite since 
$\left| \log \frac{b_{i,i+2}}{b_i} \right|  < 1$, and $b_{i,i+2}$ is finite by the 
definition of $\Phi (\ga_{i+2})$. Finally, $\lim_{j \rightarrow \infty} b_{ij} = b_i$. 

Denote $a_{ij} = b_{i,j+1} - b_{ij}$ and $a_{i0} = b_{i1}$. 
Then $\sum_{j=1}^{\infty} \sum_{i=1}^{\infty} | a_{ij} | $ is the sum of all the variations $T_i$ discussed above, 
and has been shown to be bounded by a fixed constant. 
Hence $\sum_{j=0}^{\infty} \sum_{i=1}^{\infty} | a_{ij} | =
 \Phi(\ga_1) + \sum_{j=1}^{\infty} \sum_{i=1}^{\infty} | a_{ij} |<
\infty$. Hence the sum converges absolutely. 

We have $\sum_{k=0}^{j-1} a_{ik} = b_{ij}$. Taking $j \rightarrow \infty$ yields 
$\sum_{k=0}^{\infty} a_{ik} =  \lim_{j \rightarrow \infty} b_{ij} = b_i$. Hence
$$
\Phi ( \ga ) = \sum_{i=1}^{\infty} b_i = \sum_{i=1}^{\infty} \sum_{k=0}^{\infty} a_{ik}.
$$
For any $j$, both $\Phi( \ga_j)$ and $\Phi (\ga_{j+1})$ are finite, hence 
$\Phi(\ga_{j+1})-\Phi( \ga_j)= \sum_{i=1}^{\infty} b_{i,j+1} - \sum_{i=1}^{\infty} b_{ij} = 
\sum_{i=1}^{\infty} ( b_{i,j+1} - b_{ij}) = \sum_{i=1}^{\infty} a_{ij}$.
Hence 
$$
L = \lim_{j  \rightarrow \infty} \Phi(\ga_j) = \Phi(\ga_1) +  \sum_{j=1}^{\infty} ( \Phi(\ga_{j+1})
- \Phi (\ga_j)) = \sum_{i=1}^{\infty} b_{i1} + \sum_{j=1}^{\infty} \sum_{i=1}^{\infty}  a_{ij} = 
 \sum_{j=0}^{\infty} \sum_{i=1}^{\infty}  a_{ij}.
$$ 
By the absolute convergence of $\sum_{j=0}^{\infty} \sum_{i=1}^{\infty} | a_{ij} | < \infty$,
the limit $L$ above exists, and we
can apply Fubini's theorem to obtain
$$
\Phi ( \ga ) = \sum_{i=1}^{\infty} \sum_{j=0}^{\infty} a_{ij} = 
\sum_{j=0}^{\infty} \sum_{i=1}^{\infty}  a_{ij} = L,
$$
which completes the proof.

\end{proof}

\subsection*{Finalizing the argument.} Let $\ga$ be the limit from 
the previous section. We claim the no machine works properly on 
$\ga$. For any $i$, $I_i$ is an initial segment of $\ga$ by definition. 
We only need to see that $r(\ga) \in [l_i, r_i]$ to show that $M_i^{\phi}$ 
fails on $\ga$. We claim that 
\beq
\label{rga}
 r(\ga) = \lim_{i \rightarrow \infty} r_i. 
\eeq
$r(\ga) \in [l_i, r_i]$ follows, since $l_i$ is nondecreasing 
and $r_i$ is nonincreasing. 

To show \eref{rga} observe that $\ga = \lim_{i \rightarrow \infty} \ga_i$, 
and by continuity of $f$, $f(\ga) = \lim_{i \rightarrow \infty} f(\ga_i)$.
Together with Lemma \ref{main2} we obtain
$$
\lim_{i \rightarrow \infty} \log(r_i) = 
\lim_{i \rightarrow \infty} (f(\ga_i)-\Phi(\ga_i)) = 
\lim_{i \rightarrow \infty} f(\ga_i)
-\lim_{i \rightarrow \infty} \Phi(\ga_i) = f(\ga) - \Phi(\ga) = 
\log( r(\ga)),
$$
thus \eref{rga} follows, which completes the proof. 
}

\end{document}